\documentclass{amsart}
\usepackage{latexsym,amssymb,amsmath,epsfig}
\usepackage{epsfig,graphicx,latexsym}
\usepackage{pstricks,pstricks-add,pst-math,pst-xkey}
\usepackage[latin1]{inputenc} % for Latex 2e
\def\proof{{\boldmath $Proof.$}\hskip 0.3truecm}

\newtheorem{lm}{Lemma}[section]
\def\proof{{\it Proof.}\nobreak\\}%% proof
\def\qed{\hfill$\Box$ \bigskip}%%%%% end of proof
\newtheorem{tm}{Theorem}[section]

\newtheorem{pro}{Proposition}[section]
\newtheorem{co}{Corollary}[section]

\newtheorem{rem}[tm]{Remark}

\newtheorem{ex}[tm]{Example}

\begin{document}

\title{Connectivity and other invariants of generalized products of graphs}

\author{S.C. L\'{o}pez \and F.A. Muntaner-Batle}

\maketitle

\begin{abstract}
Figueroa-Centeno et al. introduced the following product of digraphs: let $D$ be a digraph and let $\Gamma$ be a family of digraphs such that $V(F)=V$ for every $F\in \Gamma$. Consider any function $h:E(D)\longrightarrow\Gamma $.
Then the product $D\otimes_{h} \Gamma$ is the digraph with vertex set $V(D)\times V$ and
    $((a,x),(b,y))\in E(D\otimes_{h
    }\Gamma)$ if and only if $ (a,b)\in E(D)$ and $ (x,y)\in
    E(h (a,b))$.

In this paper, we introduce the undirected version of the $\otimes_h$-product, which is a generalization of the classical direct product of graphs and, motivated by it, we also recover a generalization of the classical lexicographic product of graphs that was introduced by Sabidussi en 1961. We study connectivity properties and other invariants in terms of the factors. We also present a new intersection graph that emerges when we characterize the connectivity of $\otimes_h$-product of graphs.
\end{abstract}


\section{Introduction}

We begin by introducing those concepts of classical graph theory that
will be necessary in this paper. First of all, we clarify that all the graphs considered in this
paper are assumed to be finite and, if no otherwise specified, simple. Let $G$ be a graph and let $v\in V(G)$, we let the open neighborhood of $v$ to be $N_G(v)=\{u\in V(G): uv\in E(G)\}$ and the closed neighborhood of $v$ is $N_G[v]=N_G(v)\cup \{v\}$. The degree of a vertex $v$, $|N_G(v)|$, is denoted by $d_G(v)$ and the minimum degree among the vertices of $G$ by $\delta (G)$. Let $S$ be either a subset of $V(G)$ or a subset of $E(G)$, we denote by $G[S]$ the subgraph of $G$ induced by $S$. A set $S\subset V(G)\cup E(G)$ is a {\it separating set} if its deletion, which we denote by $G-S$, disconnects $G$. The minimum size of a separating set of vertices is called the {\it connectivity} of $G$, and is denoted by $\kappa (G)$. The minimum size of a separating set of edges is called the {\it edge-connectivity} of $G$, and is denoted by $\lambda (G)$. A separating set of vertices $S$ is a $\kappa$-{\it set} if $|S|=\kappa(G)$.  Similarly, a $\lambda$-{\it set} is a separating set of edges of size $\lambda(G)$. A set $S\subset V(G)$ is a {\it dominating set} if each vertex in $V(G)\setminus S$ is adjacent to at least one vertex of $S$. A dominating set $S$ in which each vertex in $S$ has a neighbor in $S$ is called a {\it total dominating set}. The ({\it total}) {\it domination number} ($\gamma_t(G)$) $\gamma (G)$ of a graph $G$ is the minimum cardinality of a (total) dominating set. The {\it independence number} of $G$, denoted by $\alpha(G)$ is the greatest $r$ such that $rK_1$, the complement of $K_r$, is an induced subgraph of $G$. A maximal complete subgraph is a {\it clique}. The {\it clique number} $\omega (G)$ is the number of vertices of a maximum clique.
An $r$-{\it coloring} of $G$ is any function $f:V(G)\rightarrow \{0,1,2,\ldots, r-1\}$ such that if $uv\in E(G)$ then $f(u)\neq f(v)$. The {\it chromatic number} of $G$, $\chi(G)$, is the minimum $r$ for which there is an $r$-coloring of $G$.

Let $\mathcal{F}=\{S_i: \ i\in I\}$ be a family of sets. The {\it intersection graph} obtained from $\mathcal{F}$ is a graph that has a vertex $v_i$ for each $i\in I$, and for each pair $i,j\in I$, there is an edge $v_iv_j$ if and only if $S_i\cap S_j\neq \emptyset$.

Let $G$ and $H$ be two graphs. Two of the standard products of graphs are the direct and the lexicographic product. The direct product $G\otimes H$ (also denoted by $G\times H$) is the graph with vertex set $V(G)\times V(H)$ and $(a,x)(b,y)\in E(G\otimes H)$ if and only if, $ab\in E(G)$ and $xy\in E(H)$. The direct product also appears in the literature as the cross product, the categorical product, the cardinal product, the tensor product, the relational product, the Kronecker product, the weak direct product and even the cartesian product. The lexicographic product $G\circ H$ (also denoted by $G[H]$) is the graph with vertex set $V(G)\times V(H)$ and $(a,x)(b,y)\in E(G\circ H)$ if and only if, either $ab\in E(G)$ or $a=b$ and $xy\in E(H)$.

The following theorem is due to Weichsel.
\begin{tm}\cite{Wei62} \label{prod_directe_non_bipartit}
Let $G$ and $H$ be graphs with at least one edge. Then $G\otimes H$ is connected if and only if both $G$ and $H$ are connected and at least one of them is nonbipartite. Furthermore, if both are connected and bipartite, then $G\otimes H$ has exactly two connected components.
\end{tm}

%In particular, when $G\cong K_2$ we obtain the next result.
%\begin{lm}\cite{YangXu} \label{prod_directe_K2_bipartit}
%Let $H$ be a connected bipartite graph and $K_2$ be a complete graph with vertex set $\{a,b\}$. Then $K_2\otimes H$ has exactly two components. Moreover, for each $x\in V(H)$, $ax$ and $bx$ are in distinct components of $K_2\otimes H$.
%\end{lm}

\begin{rem}\label{prod_directe_K2_bipartit}
Let $a$ and $b$ be adjacent vertices of $G$.
It is not difficult to check that, if $V_1$ and $V_2$ are the stable sets of $H$ then the subgraphs induced by $(\{a\}\times V_1)\cup (\{b\}\times V_2)$ and $(\{a\}\times V_2)\cup (\{b\}\times V_1)$ are the two connected components of $G[ab]\otimes H$.
\end{rem}

The lexicographic product of two graphs $G$ and $H$, with $G$ nontrivial, is connected if and only if $G$ is connected.

Figueroa-Centeno et al. introduced the following product of digraphs in \cite{F1}: let $D$ be a digraph and let $\Gamma$ be a family of digraphs such that $V(F)=V$ for every $F\in \Gamma$. Consider any function $h:E(D)\longrightarrow\Gamma $.
Then the product $D\otimes_{h} \Gamma$ is the digraph with vertex set $V(D)\times V$ and
$((a,x),(b,y))\in E(D\otimes_{h}\Gamma)$ if and only if $ (a,b)\in E(D)$ and $ (x,y)\in E(h (a,b))$.
Notice that, when $h$ is constant, the adjacency matrix of $D\otimes_{h} \Gamma$, $A(D\otimes_{h} \Gamma)$,
coincides with the classical Kronecker product of matrices, $A(D)\otimes A(h(e))$, where $e\in E(D)$.
 When $|\Gamma|=1$, we refer to this product as the direct product of two digraphs and
 we just write $D\otimes \Gamma$ \cite{WhiRus12}. %In what follows, the direct product of two (di)graphs $D$ and $H$ will be denoted by $D\otimes H$.
The $\otimes_h$-product of digraphs has been used to establish strong relations among different labelings and specially to produce (super) edge-magic labelings for some families of graphs \cite{ILMR,LopMunRiu1,LMR_LabCons}. Some structural results can be found in \cite{AMR,LopMunRiu1,LMR_LabCons}.

An undirected version of the $\otimes_h$-product can be provided as follows: let $G$ be a graph and let $\Gamma$ be a family of graphs such that $V(F)=V$ for every $F\in \Gamma$. Consider any function $h:E(G)\longrightarrow\Gamma $.
Then the product $G\otimes_{h} \Gamma$ is the graph with vertex set $V(G)\times V$ and
$(a,x)(b,y)\in E(G\otimes_{h}\Gamma)$ if and only if $ ab\in E(G)$ and $ xy\in E(h (ab))$.

Motivated by the $\otimes_h$-product, we also recover a generalization of the lexicographic product that was introduced by Sabidussi in \cite{Sab61}. Let $G$ be a graph and let $\Gamma$ be a family of graphs. Consider any function $h:V(G)\longrightarrow\Gamma $.
Then the product $G\circ_{h} \Gamma$ is the graph with vertex set $\cup_{a\in V(G)} \{(a,x):\ x\in V(h(a))\}$ and
$(a,x)(b,y)\in E(G\circ_{h}\Gamma)$ if and only if either $ ab\in E(G)$ or $a=b$ and $xy\in E(h (a))$.

The organization of the paper is the following one. Section 2 is dedicated to connectivity of both, the $\otimes_h$-product and the generalized lexicographic product. Section 3 is focused in the study of other invariants of the generalized products, in terms of the factors. We study the independence number, the domination number, the chromatic number and the clique number. We end up this paper by presenting some structural properties in Section 4.
\section{Connectivity}
Let $G$ be a graph and let $\Gamma$ be a family of graphs such that $V(F)=V$ for every $F\in \Gamma$. Consider any function $h:E(G)\longrightarrow\Gamma $. We denote by $h(G)$ the graph with vertex set $V$ and edge set $E(h(G))=\cup_{e\in E(G)}E(h(e))$. Clearly, if $(a,x)\in V(G\otimes_h\Gamma)$ then $$d_{G\otimes_h\Gamma}(a,x)=\sum_{b\in N_G(a)}d_{h(ab)}(x).$$

%\begin{lm}
%Let $G$ be a graph and let $\Gamma =\{F_i\}_{i=1}^m$ be a family of graphs such that  $V(F_i)=V$ for every $i\in \{j\}_{j=1}^m$. Consider any function $h:E(G)\longrightarrow\Gamma $. Then $\delta (G\otimes_{h} \Gamma)=\delta (G)\min\{\delta (h(e)): \ e\in E(G)\}$.
%\end{lm}

The next result is, in some sense, a generalization of Theorem \ref{prod_directe_non_bipartit}.
\begin{tm}\label{connectivi_directe_nonbipartit}
Let $G$ be a nontrivial connected graph and let $\Gamma$ be a family of nontrivial connected graphs such that  $V(F)=V$ for every $F\in \Gamma$. Consider any function $h:E(G)\longrightarrow\Gamma $. Then, $G\otimes_{h} \Gamma$ is connected if and only if at least one of $G$ or  $h(G)$ is nonbipartite.
\end{tm}

\proof First assume that $G$ and $h(G)$ are bipartite graphs with stable sets $V(G)=A\cup B$ and $V=C\cup D$. Then, there are no edges between the sets of vertices $(A\times C)\cup (B\times D)$ and $(B\times C)\cup (A\times D)$ in
$G\otimes_{h} \Gamma$, hence $G\otimes_{h} \Gamma$ is disconnected.

Assume now that $h(G)$ is nonbipartite. Since $G$ is connected and $\delta(F)\ge 1$, for every $F\in \Gamma$, in order to prove that $G\otimes_{h} \Gamma$  is connected, we only have to prove that there exists $a\in V(G)$ such that for each pair of vertices $x,y\in V$ there is a path in $G\otimes_{h} \Gamma$ connecting $(a,x)$ and $(a,y)$. If there exists $e\in E(G)$ such that $h(e)$ is nonbipartite, then by Theorem \ref{prod_directe_non_bipartit}, the graph $G[e]\otimes h(e)$ is connected. Thus, we obtain that $G\otimes_{h} \Gamma$ is connected. Suppose now that $h(e)$ is bipartite, for each $e\in E(G)$. By Theorem \ref{prod_directe_non_bipartit}, the graph $G[e]\otimes h(e)$ has exactly two components.

Since all elements in $\Gamma$ are connected and $h(G)$ is nonbipartite, there exist $a,b,c\in V(G)$, two elements of $\Gamma$, namely $F_1, F_2$, such that $h(ab)=F_1$, $h(bc)=F_2$, and the graph $(V, E(F_1)\cup E(F_2))$ is nonbipartite. We denote by $V_1^i, V_2^i$ the stable sets of $F_i$, $i=1,2$.

{\bf\it Claim}. For each pair of vertices $x,y\in V$, there is a path connecting $(a,x)$ and $(a,y)$ in the subgraph of $G\otimes_{h} \Gamma$ induced by $\{a,b,c\}\times V$.

Suppose that $(a,x)$ and $(a,y)$ are in different components of $G[ab]\otimes h(ab)$. Otherwise, the claim is trivial. Without loss of generality suppose that $x\in V_1^1$ and $y\in V_2^1$. Assume first that $V_k^1\cap V_l^2\neq \emptyset$, for each pair $l,k\in \{1,2\}$, and let $x'\in V_1^1\cap V_2^2$ and $y'\in V_2^1\cap V_2^2$. By Remark \ref{prod_directe_K2_bipartit}, there exists a path in $G[ab]\otimes F_1$ connecting $(a,x)$ and $(b,y')$. Similarly, $(a,y)$ and $(b,x')$ are connected in $G[ab]\otimes F_1$. Moreover, Remark \ref{prod_directe_K2_bipartit} implies, since $x', y'\in  V_2^2$, that for each $z\in V_1^2$ there is a path in $G[bc]\otimes F_2$ connecting $(b,y')$ and $(c,z)$, and also a path connecting $(b,x')$ and $(c,z)$. Therefore, there is a path connecting $(a,x)$ and $(a,y)$ in the subgraph of $G\otimes_{h} \Gamma$ induced by $\{a,b,c\}\times V$. Assume now that $V_k^1\cap V_l^2= \emptyset$, for some pair $l,k\in \{1,2\}$. Without loss of generality assume that $V_1^1\cap V_1^2=\emptyset$. Then, we have $V_2^1\cap V_2^2\neq \emptyset$. Otherwise, the graph $(V, E(F_1)\cup E(F_2))$ is bipartite, a contradiction. We then proceed as in the above case.

Assume now that $h(G)$ is bipartite and that $G$ is nonbipartite. Let $C=a_0a_1\ldots a_{2k}a_0$ be an odd cycle in $G$. Since $h(G)$ is bipartite, it follows that there exists a partition $V=V_1\cup V_2$, such that $V_1$ and $V_2$ are the stable sets of $h(a_ia_{i+1})$ and $h(a_{2k}a_0)$, for each $i=0,1,\ldots, 2k-1$. Which implies, since the cycle is odd, that we can connect every vertex in $\{a_0\}\times V_1$ to every vertex in $\{a_1\}\times V_2$, $\{a_2\}\times V_1$, and so on, until, $\{a_{2k}\}\times V_1$ and finally, $\{a_0\}\times V_2$. Therefore and by Remark \ref{prod_directe_K2_bipartit}, the graph $G\otimes_h \Gamma$ is connected.
\qed

When we do not assume that the elements of $\Gamma$ are connected we can find examples of disconnected graphs that are of the form $G\otimes_h\Gamma$, with both $G$ and $h(G)$ nonbipartite and connected.

\begin{ex}
Let $V=\{x,y,z,t\}$. Consider the graphs $F_i$ on $V$, $i=1,2,3$, defined by, $E(F_1)=\{xz,yt\}$, $E(F_2)=\{xy,zt\}$ and $E(F_3)=\{xt,yz\}$. Let $h: E(C_3)\rightarrow \{F_i\}_{i=1}^3$ be any bijective function. Then, $C_3\otimes_h\{F_i\}_{i=1}^3\cong 4C_3$. However, both graphs, $C_3$ and $h(C_3)\cong K_4$, are nonbipartite and connected.
\end{ex}

The next results give sufficient conditions to guarantee connectivity in $G\otimes_h \Gamma$ when the family $\Gamma $ contains disconnected graphs. The first result appears in the proof of Theorem \ref{connectivi_directe_nonbipartit}. Notice that, if there exists $ab\in V(G)$ such that $h(ab)$ has an isolated vertex and either $a$ or $b$ is a vertex of degree $1$ in $G$, then the graph $G\otimes_h\Gamma$ is not connected. So, in what follows, we assume that all vertices of $F$ have degree at least $1$, for every $F\in \Gamma$. Recall that, since $G$ is connected and $\delta(F)\ge 1$, for every $F\in \Gamma$, in order to prove that $G\otimes_{h} \Gamma$ is connected, we only have to prove that there exists $a\in V(G)$ such that for every pair $x,y\in V$ there is a path in $G\otimes_{h} \Gamma$ connecting $(a,x)$ and $(a,y)$. This fact is guaranteed in the following two lemmas.

\begin{lm}\label{one_nonbipartit_connected}
Let $G$ be a nontrivial connected graph and let $\Gamma$ be a family of graphs such that  $V(F)=V$ and $\delta(F)\ge 1$, for every $F\in \Gamma$. Consider any function $h:E(G)\longrightarrow\Gamma $. If there exists $e\in E(G)$ such that $h(e)$ is nonbipartite and connected.
Then, the graph $G\otimes_{h} \Gamma$ is connected.
\end{lm}

\proof Let $e\in E(G)$ such that $h(e)$ is nonbipartite and connected. Then, by Theorem \ref{prod_directe_non_bipartit}, the graph $G[e]\otimes h(e)$ is connected. \qed

\begin{lm}\label{one_bipartit_connected}
Let $G$ be a nontrivial connected graph and let $\Gamma$ be a family of graphs such that  $V(F)=V$ and $\delta(F)\ge 1$, for every $F\in \Gamma$. Consider any function $h:E(G)\longrightarrow\Gamma $. Let $ab, bc\in E(G)$ such that $h(ab)$ is bipartite and connected, with stable sets $V_1$ and $V_2$ and assume that at least one of the following holds:

\begin{enumerate}
  \item[(i)] One of the components of $h(bc)$ is nonbipartite and contains vertices of $V_1$ and $V_2$.
  \item[(ii)] One of the components of $h(bc)$ is bipartite, but one of the stable sets contains vertices of $V_1$ and $V_2$.
\end{enumerate}

Then, the graph $G\otimes_{h} \Gamma$ is connected.
\end{lm}
\proof By Theorem \ref{prod_directe_non_bipartit}, the subgraph $G[ab]\otimes h(ab)$ has two components, which are the subgraphs induced by the sets of vertices, $(\{a\}\times V_1)\cup (\{b\}\times V_2)$ and $(\{a\}\times V_2)\cup (\{b\}\times V_1)$. Let us prove that (i) implies that $G\otimes_h\Gamma$ is connected. Let $C_{bc}$ be a nonbipartite component of $h(bc)$ and $x,y\in V(C_{bc})$ with $x\in V_1$ and $y\in V_2$. Since $C_{bc}$ is nonbipartite, the subgraph of $G\otimes_h\Gamma$ induced by $\{b,c\}\times V(C_{bc})$ is connected and contains vertices of the two components of $G[ab]\otimes h(ab)$. Therefore, all vertices of the form $\{a,b\}\times V$ are in the same component of $G\otimes_h\Gamma$ and the result follows. (ii) Suppose now that $C_{bc}$ is a bipartite component of $h(bc)$ which contains two vertices $x,y\in V(h(bc))$ in the same stable set, namely $V_1(C_{bc})$, such that $x\in V_1$ and $y\in V_2$. Thus, the subgraph of $G\otimes_h\Gamma$ induced by $\{b\}\times V_1(C_{bc})$ connects vertices of $(\{a\}\times V_1)\cup (\{b\}\times V_2)$ with vertices of $(\{a\}\times V_2)\cup (\{b\}\times V_1)$. Hence, the subgraph induced by $\{a,b\}\otimes V$ belongs to the same component of $G\otimes_h\Gamma$. Therefore, the result is proved.
\qed

The next result is a technical theorem that presents an interesting relation between some properties of partitions and the connectivity of the intersection graph obtained from them.

\begin{tm}\label{cicle_partitions}
Let $\mathcal{P}_1(A)$, $\mathcal{P}_2(A)$, \ldots, $\mathcal{P}_m(A)$  be partitions of a set $A$ and $G$ the intersection graph obtained from $\cup_{i=1}^m\mathcal{P}_i(A)$. Then, $G$ is disconnected
 if and only if, there exists nonempty $\mathcal{A}_i\subset \mathcal{P}_i(A)$, for each $i\in \{1,2,\ldots, m\}$, such that  $\cup_{A_i\in \mathcal{A}_i}A_i= \cup_{A_j\in \mathcal{A}_j}A_j\neq V$, for each $i,j$ with $1\le i\le j\le m$.
\end{tm}

\proof We denote each vertex of $G$ with the name of the corresponding set of $\mathcal{P}_1(A)\cup \mathcal{P}_2(A)\cup\ldots \cup\mathcal{P}_m(A)$. Let us see the sufficiency. Assume that, for each $k$ with $1\le k\le m$, there exists $\mathcal{A}_k\subset \mathcal{P}_k(A)$, such that  $\cup_{A_i\in \mathcal{A}_i}A_i= \cup_{A_j\in \mathcal{A}_j}A_j$, for each $i,j$ with $1\le i\le j\le m$. Let $A_i\in \mathcal{A}_i$ and $B_j\in \mathcal{P}_j(A)\setminus \mathcal{A}_j$, for each $i,j$ with $1\le i\le j\le m$. By hypothesis, we have that $A_i\cap B_j=\emptyset$ and thus, $A_iB_j\notin E(G)$. Hence, the subgraphs $G[\cup_{i=1}^m \mathcal{A}_i]$ and $G[\cup_{i=1}^m (\mathcal{P}_i(A)\setminus \mathcal{A}_i)]$ are not in the same connected component of $G$.

Let us prove now the necessity. Suppose that $G$ is disconnected, $H$ is a connected component of $G$ and let $V(H)\cap  \mathcal{P}_i(A)=\mathcal{A}_i$, for $i=1,2,\ldots,m$. Clearly, we have that $\mathcal{A}_i\neq \mathcal{P}_i(A)$, otherwise $G$ is connected. We will prove that $\cup_{A_i\in \mathcal{A}_i}A_i= \cup_{A_j\in \mathcal{A}_j}A_j$. We proceed by contradiction. Assume to the contrary that $\cup_{A_i\in \mathcal{A}_i}A_i\neq\cup_{A_j\in \mathcal{A}_j}A_j$, for some pair $i,j$ with $1\le i\le j\le m$. Without loss of generality assume that there exists $a\in A_i$ such that $a\notin \cup_{A_j\in \mathcal{A}_j}A_j$. Thus, since we are dealing with partitions of $A$, there exists $B_j\in  \mathcal{P}_j(A)\setminus \mathcal{A}_j$ such that $a\in B_j$, and hence, $A_iB_j\in E(H)$. Therefore, we have that $B_j\in V(H)$, contradicting the fact that $V(H)\cap  \mathcal{P}_j(A)=\mathcal{A}_j$.\qed

For every graph $F$, we can associate a partition of $V(F)$, namely $\mathcal{P}_F(V(F))$ as follows. If $H$ is a bipartite component of $F$, then each stable set of $V(H)$ is an element of $\mathcal{P}_F(V(F))$. Otherwise, the set $V(H)$ itself is an element of $\mathcal{P}_F(V(F))$.

Next, we are ready to state the following results. The proofs are a direct consequence of Theorem \ref{cicle_partitions}. The first result gives a characterization for the connectivity of $G\otimes_h\Gamma$ when we relax the condition on connectedness over all elements of the family $\Gamma$, under the assumption that all connected components are nonbipartite.

\begin{tm}
Let $G$ be a nontrivial connected graph and let $\Gamma$ be a family of graphs such that  $V(F)=V$, for every $F\in \Gamma$. Consider any function $h:E(G)\longrightarrow\Gamma $. Assume that all connected components of $h(e)$ are nonbipartite, for every $e\in E(G)$. Then $G\otimes_h\Gamma$ is disconnected, if and only if, for each $e\in E(G)$ there exists $\mathcal{A}_{e}
\subset \mathcal{P}_{h(e)}(V)$, such that
 $\cup_{A_{e}\in\mathcal{A}_{e}}A_{e}=\cup_{A_{e'}\in\mathcal{A}_{e'}}A_{e'}\neq V,$
 for every $e,e'\in E(G)$.

\end{tm}

If we concentrate on the star then we can obtain a complete characterization, which does not depend on the bipartitness of the elements of $\Gamma$.

\begin{tm}
Let $G\cong K_{1,n}$ and let $\Gamma=\{F_i\}_{i=1}^m$ be a family of graphs such that $V(F_i)=V$, for each $1\le i\le m$. Assume that $h:V(H)\rightarrow \Gamma$ is a surjective mapping. Then the graph $G\otimes_h\Gamma$ is disconnected, if and only if, there exists $\mathcal{A}_{i}
\subset \mathcal{P}_{F_i}(V)$, for each $i\in \{1,2,\ldots, m\}$, such that
 $\cup_{A_{i}\in\mathcal{A}_{i}}A_{i}=\cup_{A_{j}\in\mathcal{A}_{j}}A_{j}\neq V,$
 for every $i$ and $j$ with $1\le i\le j\le m$.
\end{tm}

The above results give sufficient conditions that guarantee the connectivity of $G\otimes_h\Gamma$, when $\Gamma$ contains disconnected graphs. In order to provide a characterization, we introduce an intersection graph that we obtain from $G$, $\Gamma$ and the function $h: E(G)\rightarrow \Gamma$. Let $G$ be a nontrivial connected graph and let $\Gamma$ be a family of graphs such that $V(F)=V$, for every $F\in \Gamma$. Consider any function $h:E(G)\longrightarrow\Gamma $. For each bipartite component $C$ of a $h(ab)$, we let $$S_a(C)=(\{a\}\times V_1(C))\cup (\{b\}\times V_2(C))\ \hbox{and} \ S_b(C)=(\{a\}\times V_2(C))\cup (\{b\}\times V_1(C)),$$ where $V_1(C)$ and $V_2(C)$ are the stable sets of $C$.
If $C$ is a nonbipartite component of $h(ab)$, then we let $S_a(C)=S_b(C)= \{a,b\}\times V(C).$ We denote by $\mathcal{F}(G,\Gamma,h)$ the family of sets of $V(G)\times V$ defined by: $$\mathcal{F}(G,\Gamma,h)=\bigcup_{\tiny\begin{array}{l}
                                              a\in V(G) \\
                                              b\in N_G(a)
                                            \end{array}}\{S_{a}(C):\ C \ \hbox{is a component of}\ h(ab)\}.$$
Thus, by definition, we obtain the next characterization.

\begin{tm}\label{th_charac_conne_generalized_otimes_product}
Let $G$ be a nontrivial connected graph and let $\Gamma$ be a family of graphs such that $V(F)=V$, for every $F\in \Gamma$. Consider any function $h:E(G)\longrightarrow\Gamma $.
Then, $G\otimes_h\Gamma$ is connected if and only if the intersection graph obtained from $\mathcal{F}(G,\Gamma,h)$ is connected.
\end{tm}

\begin{rem}
It is worthly to mention that, as it happens with the direct product, the most natural setting for the $\otimes_h$-product is the class of graphs with loops. If we accept to apply our definition to graphs with loops, a rewiew of the proofs reveals that all results presented in this section remain valid in the class of graphs with loops.
\end{rem}

\subsection{Connectivity in $G\circ_h\Gamma$}
Let $G$ be a graph and let $\Gamma$ be a family of graphs. Consider any function $h:V(G)\longrightarrow\Gamma $. Clearly, if $(a,x)\in V(G\circ_h\Gamma)$ then $$\displaystyle d_{G\circ_h\Gamma}(a,x)=\sum_{b\in N_G(a)}|V(h(b))|+d_{h(a)}(x).$$
In particular, if $V(F)=V$ for every $F\in \Gamma$ then

\begin{equation}\label{min_degre_gen_lexico}
\delta (G\circ_h\Gamma)=\delta (G)|V|+\min_{\tiny\begin{array}{c}
v\in V(G) \\
d_G(v)=\delta(G)
\end{array}
}\delta(h(v)).
\end{equation}

The next result is trivial.
\begin{lm}
Let $G$ be a nontrivial graph and let $\Gamma$ be a family of graphs. Consider any function $h:V(G)\longrightarrow\Gamma $. Then $G\circ_h\Gamma$ is connected if and only if $G$ is connected.
\end{lm}

When $V(F)=V$ for every $F\in \Gamma$, we can obtain exact formulas for the connectivity and the edge-connectivity of $G\circ_h\Gamma$, which generalize the ones corresponding to $G\circ H$. The proofs are similar to the case $G\circ H$ in \cite{YangXu}. For each $a\in V(G)$, the $V$-fiber of $G\circ_h \Gamma$ with respect to $a$, refers to $_aV=\{(a,x): x\in V\}$.

\begin{tm}\label{gen_lexi_connectivity}
Let $G$ be a connected graph of order $n$ and let $\Gamma$ be a family of graphs such that $V(F)=V$, for every $F\in \Gamma$. Consider any function $h:V(G)\longrightarrow\Gamma $. Then
$$\kappa( G\circ_h\Gamma)=\left\{\begin{array}{ll}
                               (n-1)|V|+\min_{v\in V(G)}\kappa(h(v)), & G= K_n, n\ge 1,\\
                              \kappa (G)|V|, & G\neq K_n.
                            \end{array}\right.
$$
\end{tm}

\proof Suppose first that $G= K_n$, for $n\ge 1$. We claim that $\kappa (G\circ_h\Gamma)=(n-1)|V|+\min_{v\in V(G)}\kappa(h(v))$. Let $S$ be a separating set of vertices of $G\circ_h\Gamma$. If $(a,x)\notin S$ then $V(G\circ_h\Gamma)\setminus S\subset\ _aV$, otherwise, since $G$ is the complete graph, we obtain that $G\circ_h\Gamma-S$ is connected, a contradiction. Let $S_a=S\cap\ _aV$. We have that $\pi_2(S_a)$ is a separating set of $h(a)$, where $\pi_2: V(G)\times V\rightarrow V$ is the natural projection defined by $\pi_2(a,x)=x$. Therefore, $|S|\ge (n-1)|V|+\kappa(h(a))$. Let $a\in V(G)$ such that $\kappa (h(a))=\min_{v\in V(G)}\kappa(h(v))$ and consider a $\kappa$-set $S'$ of $h(a)$. Then, $((V(G)\setminus \{a\})\times V)\cup (\{a\}\times S')$ is a separating set of $G\circ_h\Gamma$. Hence, we obtain that $\kappa(G\circ_h\Gamma)\le(n-1)|V|+\min_{v\in V(G)}\kappa(h(v))$ and the claim follows.

Suppose now that $G$ is a connected graph of order $n$ not isomorphic to $K_n$, for $n\ge 1$. Let $S$ be a $\kappa$-set of $G$. Clearly, $S\times V$ is a separating set of $G\circ_h\Gamma$ and thus, $\kappa(G\circ_h\Gamma)\le \kappa(G)|V|$. Consider now a separating set of vertices $S$ of $G\circ_h\Gamma$, we will show that $|S|\ge \kappa(G)|V|$. We claim that there exist two vertices $(a,x)$ and $(b,y)$ that are in different connected components of $G\circ_h\Gamma- S$, with $a\neq b$. Suppose to the contrary that $G\setminus S\subset \ _aV$, for some $a\in V(G)$. We obtain that $|S|\ge (n-1)|V|$, a contradiction with the fact that $\kappa(G\circ_h\Gamma)\le \kappa(G)|V|$ and $G\neq K_n$. Since $a\neq b$ there exist $\kappa(G)$ disjoint paths in $G$, namely, $P_1, P_2,\ldots, P_{\kappa(G)}$ connecting $a$ and $b$. Let $P=aa_1a_2\ldots a_rb$ be one of these paths. If for all $a_i$ there exists $x_i\in V$ such that $(a_i,x_i)\notin S$, then $(a,x)$ and $(b,y)$ are connected through $(a_i,x_i)$, for $i=1,2,\ldots, r$. Thus, for each path $P_j$, there exists $a_i$, such that $\{a_i\}\times V\subset S$. Hence, we have that $|S|\ge \kappa(G)|V|$ and the equality $\kappa(G\circ_h\Gamma)=|\kappa(G)||V|$ is proved when $G\neq K_n$.\qed

Notice that, the previous proof also shows that, when $G\ne K_n$, a $\kappa$-set of $G\circ_h\Gamma$ is of the form $\cup_{a\in S}\ _aV$, where $S$ is $\kappa$-set of $G$. If we remove the hypothesis $V(F)=V$, for every $F\in \Gamma$ then, we cannot obtain this conclusion. However, the following inequality still holds: $$\kappa(G\circ_h\Gamma)\le \min_{S}\sum_{v\in S}|V(h(v)|,$$ where the minimum is taken over all separating sets of vertices of $G$.

\begin{tm}\label{gen_lexi_edge_connectivity}
Let $G$ be a connected graph of order $n\ge 2$ and let $\Gamma$ be a family of nontrivial graphs such that $V(F)=V$, for every $F\in \Gamma$. Consider any function $h:V(G)\longrightarrow\Gamma $. Then
$$\lambda( G\circ_h\Gamma)=\min\{\lambda (G)|V|^2, \delta(G\circ_h\Gamma)\}.
$$
\end{tm}

\proof Clearly, all edges incident to a vertex (of minimum degree) form a separating set of edges. Similarly, from a $\lambda$-set $S$ of $G$, we obtain a separating set $\{(a,x)(b,y):\ ab\in S\ \hbox{and}\ x,y\in V\}$ of $G\circ_h\Gamma$ of size $|S||V|^2$. Thus, we have that $\lambda( G\circ_h\Gamma)\le\min\{\lambda (G)|V|^2, \delta(G\circ_h\Gamma)\}$.

Let $S$ be a $\lambda$-set of $G\circ_h\Gamma$. Then, $G\circ_h\Gamma-S$ has exactly two connected components, namely $C_1$ and $C_2$. Consider the subsets $A=\{a\in V(G):\ _aV\cap V(C_1)\neq\emptyset\}$ and $B=\{b\in V(G):\ _bV\cap V(C_2)\neq\emptyset\}$. We can assume that $A\cap B\neq \emptyset$, otherwise, $A\cup B$ is a partition of $V(G)$ and, thus, the cardinality of the set of edges joining vertices of $A$ with vertices of $B$ is at least, $\lambda (G)$. Hence, we obtain that $|S|\ge \lambda (G)|V|^2$, and the result follows.

Let $a\in A\cap B$. For every $b\in N_G(a)$, denote by $G_b[a,b]$ the bipartite subgraph of $G\circ_h\Gamma$ with stable sets $_aV$ and $_bV$. By definition of $\circ_h$, the graph $G_b[a,b]$ is a bipartite complete graph, and thus, with edge connectivity $|V|$.
Let $S_{ab}=S\cap E(G_b[a,b])$ and $S_a=S\cap E(G\circ_h\Gamma[_aV])$, where $G\circ_h\Gamma[_aV]$ is the subgraph of $G\circ_h\Gamma$ induced by $_aV$. Since the graph $G_b[a,b]$ has connectivity $|V|$, we have that
\begin{equation}\label{ineq_1}
    |S_{ab}|\ge |V|, \ \hbox{for all} \ b\in N_G(a).
\end{equation}
Moreover, we claim that,
\begin{eqnarray}\label{ineq_deltapartition}
  |S_{ab}|+|S_a| &\ge & \delta (h(a))+|V|.
\end{eqnarray}
In order to prove inequality (\ref{ineq_deltapartition}), we consider the sets $X=\ _aV\cap V(C_1)$ and $Y=\ _aV\cap V(C_2)$. Without lost of restriction assume that $|X|\le |Y|$. We can suppose that $|X|\ge 2$, otherwise, $|S_a|\ge \delta (h(a))$ and inequality (\ref{ineq_deltapartition}) holds. Let $X=\{(a,x_1),(a,x_2),\ldots, (a,x_r)\}$ and $\{(a,y_1),(a,y_2),\ldots,(a,y_r)\}\subset Y$. Then, there are $|V|$ edge-disjoint paths in $G_b[a,b]$ joining $(a,x_i)$ and $(a,y_i)$, for every $i=1,2,\ldots, r$. Thus, $|S_{ab}|\ge |X||V|\ge 2|V|>\delta (h(a))+|V|$, and inequality (\ref{ineq_deltapartition}) is proved. Hence, using that
$|S|\ge |S_a|+\sum_{b\in N_G(a)}|S_{ab}|=|S_a|+|S_{ab}|+\sum_{b'\in N_G(a)\setminus\{b\}}|S_{ab'}|$ and inequalities (\ref{ineq_1}) and (\ref{ineq_deltapartition}), we obtain that either $|S|\ge \delta (h(a))+\delta(G)|V|$, when $d_G(a)=\delta_G$, or $|S|\ge (\delta (G)+1)|V|$, otherwise. Therefore, using (\ref{min_degre_gen_lexico}), we have that $|S|\ge \delta(G\circ_h\Gamma)$ and the result follows.
\qed

%Note that, the previous proof also shows that we can compute the exact value of $\lambda(G\circ_h \Gamma)$ if we assume that, for all $x\in G$ there exists $y\in N(x)$ such that $|V(h(y))|\ge \delta (h(x))$.
%
%\begin{tm}\label{gen_lexi_edge_connectivity}
%Let $G$ be a connected graph of order $n\ge 2$ and let $\Gamma$ be a family of nontrivial graphs. Consider any function $h:V(G)\longrightarrow\Gamma $ such that,  for all $x\in G$ there exists $y\in N(x)$ such that $|V(h(y))|\ge \delta (h(x))$. Then
%$$\delta (G\circ_h\Gamma)=\min\{\delta (G\circ_h\Gamma),\min_{S}\sum_{xy\in S} |V(h(x))||V(h(y))|\},$$
%where the second minimum is taken over all $\lambda$-sets of $G$.
%\end{tm}

\section{Other invariants of generalized products}
In this section we study some invariants related to the generalized products $\otimes_h$ and $\circ_h$. We start with the independence number. Based on Proposition 8.10 in \cite{ImrKlav00}, which is related to the independence number of the direct product, we have a clear lower bound for the independence number of $G\otimes_h \Gamma$. For each $a\in V(G)$, the $V$-fiber of $G\otimes_h \Gamma$ with respect to $a$, refers to $_aV=\{(a,x): x\in V\}$ and the $G$-fiber of $G\otimes_h \Gamma$ with respect to $x\in V$ is $G_x=\{(a,x): a\in V(G)\}$.

\begin{pro}
Let $G$ be a graph and let $\Gamma$ be a family of graphs such that $V(F)=V$, for every $F\in \Gamma$. Consider any function $h:E(G)\longrightarrow\Gamma $. Then,
$\alpha(G\otimes_h \Gamma)\ge \max\{\alpha (G)|V|,\alpha (h(G))|V(G)|\}.$
\end{pro}

\proof The inclusion $E(G\otimes_h\Gamma)\subset E(G\otimes h(G))$ implies that $\alpha(G\otimes_h \Gamma)\ge \alpha (G\otimes h(G))$. Suppose now that $I$ is an independent set of $G$, then $\cup_{a\in I}\ _aV$ is an independent set of $G\otimes h(G)$. Similarly, if $J$ is an independent set of $h(G)$ then  $\cup_{x\in J} G_x$ is an independent set of $G\otimes h(G)$. Therefore, we get the result.\qed

With respect to the independence number in $G\circ_h\Gamma$ we can obtain an exact formula in terms of the independent sets of $G$.

\begin{pro}\label{gen_lexi_independence_number}
Let $G$ be a graph of order $n\ge 2$ and let $\Gamma$ be a family of graphs. Consider any function $h:V(G)\longrightarrow\Gamma $. Then
$$\alpha (G\circ_h\Gamma)=\max_{S}\sum_{a\in S} \alpha(h(a)),$$
where the maximum is taken over all independent sets of vertices $S$ of $G$.
\end{pro}

\proof Let $S$ be an independent set of vertices of $G$ and let $I_a$ be a set of independent vertices of $h(a)$. Then, the disjoint union $\cup_{a\in S}\{(a,x): x\in I_a\}$ is an independent set of $G\circ_h\Gamma$. Thus, we have that $\alpha(G\circ_h\Gamma)\ge \max_{S}\sum_{a\in S} \alpha(h(a))$.
Suppose now that $S^\circ$ is a maximal independent set of vertices of $G\circ_h\Gamma$ and let $S_a=\{x\in V(h(a)): (a,x)\in S^\circ\}$. For any pair $a,b$ of vertices of $G$ such that $S_a$ and $S_b$ are nonempty, we have that $a$ and $b$ are independent vertices. Moreover, the maximality of $S^\circ$ implies that if $|S_a|\ge 1$ then $|S_a|=|\alpha(h(a))|$. Hence, we obtain that
$|S^\circ|= \sum_{a\in S} \alpha(h(a))$, where $S$ is some independent set of vertices of $V(G)$. Therefore, we have that $\alpha(G\circ_h\Gamma)\le \max_{S}\sum_{a\in S} \alpha(h(a))$ and the result is proved.
\qed

\subsection{Domination number}

%Let $G$ be a graph. For a total dominating set $D$ of $G$ let $\gamma_D(G)$ denote the size of the smallest subset of $D$ that dominates $G$.

Although Gravier and Khellady \cite{GraKhe95} posed a kind of Vizing's conjecture for the direct product of graphs, namely $\gamma (G\otimes H)\ge \gamma (G)\gamma (H)$, a year later Nowakowski and Rall \cite{NowRall96} gave a counterexample. In fact, Klav$\check{z}$ar and Zmazek \cite{KlaZam} showed that the difference $\gamma (G)\gamma (H)-\gamma (G\otimes H)$ can be arbitrarily large. Recently, Meki$\check{s}$ has shown in \cite{Meki10} that for arbitrary graphs $G$ and $H$, we have $\gamma (G\otimes H)\ge \gamma (G)+\gamma (H)-1$. Thus, since $E(G\otimes_h\Gamma)\subset E(G\otimes h(G))$, we obtain the next easy corollary. Recall that $h(G)$ is the graph with vertex set $V$ and edge set $\cup_{e\in E(G)}E(h(e))$.

\begin{co}
Let $G$ be a graph and let $\Gamma$ be a family of graphs such that $V(F)=V$, for every $F\in \Gamma$. Consider any function $h:E(G)\longrightarrow\Gamma $. Then,
$\gamma (G\otimes_h \Gamma)\ge \gamma (G)+\gamma (h(G))-1.$
\end{co}

Inspirated by Meki$\check{s}$' lower bound proof, we improve the above lower bound for the domination number of the $\otimes_h$-product. We let $ h(G^a)=(V, \cup_{b\in N_{G}(a)} E(h(ab)))$, for each $a\in V(G)$.
\begin{tm}
Let $G$ be a graph and let $\Gamma$ be a family of graphs such that $V(F)=V$, for every $F\in \Gamma$. Consider any function $h:E(G)\longrightarrow\Gamma $. Then,
$\gamma (G\otimes_h\Gamma)\ge \gamma (G)+\min_{a\in V(G)}\gamma (h(G^a))-1.$

\end{tm}

\proof Notice that, if $D\subset V(G)\times V$ is a dominating set of $G\otimes_h\Gamma$, then $\pi_1(D)$ is a dominating set of $G$, where $\pi_1:V(G)\times V\rightarrow V(G)$ defined by $\pi_1(a,x)=a$. In particular, $\gamma(G)\le |\pi_1(D)|$. Similarly, for each $a\in V(G)$, we can check that the set $\pi_2 (\cup_{b\in N_G[a]}D\cap\ _bV)$ is a dominating set of $h(G^a)$, where $\pi_2$ is as defined in the proof of Theorem \ref{gen_lexi_connectivity}. Indeed, for every $x\in V\setminus \pi_2 (\cup_{b\in N_G[a]}D\cap\ _bV)$, the vertex $(a,x)$ is adjacent to some $(b,y)\in D$. Thus, $ab\in E(G)$ and $xy\in E(h(ab))$ and we obtain that $\gamma (h(G^a))\le |\pi_2(\cup_{b\in N_G[a]}D\cap\ _bV)|$.

Assume to the contrary that there exists a dominating set $D\subset V(G)\times V$ of $G\otimes_h\Gamma$ with $|D|=\gamma (G)+\min_{a\in V(G)}\gamma (h(G^a))-2$. If $\gamma (G)=1$ then $|D|=\gamma(h(G^a))-1$ and $|\pi_2(\cup_{b\in N_G[a]}D\cap\ _bV)|\le \gamma(h(G^a)-1$, a contradiction. Similarly, if $\min_{a\in V(G)}\gamma(h(G^a))=1$, the set $\pi_1(D)$ gives a dominating set of $G$ of size at most $\gamma (G)-1$, also a contradiction. Thus, we may assume that $\gamma(G), \gamma (h(G^a))\ge 2$, for each $a\in V(G)$. Let $D_0=\{(a_1,x_1),(a_2,x_2),\ldots, (a_{\gamma(G)-1},x_{\gamma(G)-1})\}$ be a proper subset of $D$, with $a_i\ne a_j$, for each pair $i,j$, with $i\neq j$. Since $|\pi_1(D_0)|=\gamma(G)-1$, there exists $c\in V(G)\setminus \pi_1(D_0)$ which is not adjacent to any of the vertices of $\pi_1(D_0)$.
Consider now the set $D\setminus D_0$. Since $|D\setminus D_0|=\min_{a\in V(G)}\gamma(h(G^a))-1$, we have that $|\pi_2(\cup_{b\in N_G[c]}(D\setminus D_0)\cap\ _bV)|\le \gamma(h(G^c))-1$. Thus, there exists $y\in V\setminus \pi_2(\cup_{b\in N_G[c]}(D\setminus D_0)\cap\ _bV)$ which is not adjacent to any of the vertices of $\pi_2(\cup_{b\in N_G[c]}(D\setminus D_0)\cap\ _bV)$. Notice that the vertex $(c,y)$ is not adjacent to any vertex of $D$, which implies, since $D$ is a dominating set that $(c,y)\in D$. Moreover, the condition $c\notin \pi_1(D_0)$ implies that, $c\in D\setminus D_0$ and that $y\in\pi_2(\cup_{b\in N_G[c]}(D\setminus D_0)\cap\ _bV)$, a contradiction.\qed

In particular, since for each graph $G$, we have $\gamma_t(G)\ge \gamma (G)$, we also obtain the following result.
\begin{co}
Let $G$ be a graph and let $\Gamma$ be a family of graphs such that $V(F)=V$, for every $F\in \Gamma$. Consider any function $h:E(G)\longrightarrow\Gamma $. Then,
$\gamma_t (G\otimes_h\Gamma)\ge \gamma (G)+\min_{a\in V(G)}\gamma (h(G^a))-1.$

\end{co}

With respect the upper bound for the domination number of direct product of graphs, Bre$\check{s}$ar et al. proved in \cite{BreKlaRal} the next theorem.

\begin{tm}\cite{BreKlaRal}\label{theo_dominating_direct_product}
Let $G$ and $H$ be arbitrary graphs. Then
$\gamma (G\otimes H)\le 3\gamma(G)\gamma(H).$
\end{tm}

The next generalization can be trivially obtained from Theorem \ref{theo_dominating_direct_product}.

\begin{co}
Let $G$ and $F$ be graphs and let $\Gamma$ be a family of graphs such that $V(F')=V(F)=V$ and $F$ is a subgraph of $F'$, for every $F'\in \Gamma$. Consider any function $h:E(G)\longrightarrow\Gamma $. Then,
$\gamma(G\otimes_h \Gamma)\le 3\gamma (G)\gamma(F).$
\end{co}
\proof Let us consider the inclusion $E(G\otimes F)\subset E(G\otimes_h\Gamma)$. Thus, we have that $\gamma(G\otimes_h\Gamma)\le \gamma(G\otimes F)$ and the result follows from Theorem \ref{theo_dominating_direct_product}.\qed

If the existence of a spanning connected subgraph $F$ for every graph $F'$ in $\Gamma$ is not assumed in $\Gamma$, then, we cannot control the size of a dominating set of $G\otimes_h \Gamma$. However, a similar proof as the one of Theorem \ref{theo_dominating_direct_product} in \cite{BreKlaRal} allows us to obtain the next result.

\begin{lm}
Let $G$ be a graph and let $\Gamma$ be a family of graphs
such that $V(F)=V$, for every $F\in \Gamma$. Consider any function $h:E(G)\longrightarrow\Gamma $. Assume that $D$ is a total dominating set of $G$ and let $D_e$ be a total dominating set of $h(e)$. Let $A\subset D$ be a dominating set of $G$ and let $B_e\subset D_e$ be a dominating set of $h(e)$.
Then $X=(A\times \cup_eD_e)\cup (D\times \cup_eB_e)$ is a dominating set of $G\otimes_h \Gamma$. Thus, $\gamma(G\circ_h\Gamma)\le |X|$.
\end{lm}
\proof We let $D'=\cup_{e\in E(G)}D_e$ and $B=\cup_{e\in E(G)}B_e$. We claim that $X=(A\times D')\cup (D\times B)$ is a dominating set of $G\otimes_h \Gamma$. We consider three cases.
First, assume that $a\in D\setminus A$ and $x\in D'\setminus B$. Since $A$ is a dominating set of $G$ there exists $b\in A$ such that $ab\in E(G)$. In particular, we have that  $x\notin B_{ab}$. Thus, there is $y\in B_{ab}$ such that $xy\in E(h(ab))$. Hence, we obtain that $(b,y)\in X$ and $(a,x)(b,y)\in E(G\otimes_h \Gamma)$.
Assume now that $a\in V(G)\setminus D$ and $x\in V$. Since $A$ dominates $G$ there is $b\in A$ such that $ab\in E(G)$. Now, the existence of $y\in D'$ such that $xy\in E(h(ab))$ is guaranteed by considering the total dominating set $D_{ab}$ of $h(ab)$.
Finally, assume that $(a,x)\in V(G)\times (V\setminus D')$. Since $D$ is a total dominating set there exists $b\in D$ such that $ab\in E(G)$. In particular, we have that  $x\notin B_{ab}$. Thus, there is $y\in B_{ab}$ such that $xy\in E(h(ab))$. Hence, we obtain that $(b,y)\in X$ and $(a,x)(b,y)\in E(G\otimes_h \Gamma)$.\qed

\subsubsection{Domination number in $G\circ_h\Gamma$}

Nowakowski and Rall proved in \cite{NowRall96} the inequality $\gamma(G\circ H)\le \gamma (G)\gamma (H)$. This inequality can be generalized to the $\circ_h$-product as follows

\begin{lm}
Let $G$ be a graph and let $\Gamma$ be a family of graphs. Consider any function $h:V(G)\longrightarrow\Gamma $. Then
$\gamma (G\circ_h\Gamma)\le \min_D \sum_{a\in D}\gamma (h(a)),$ where the minimum is taken over all dominating sets $D$ of $G$.
\end{lm}

\subsection{The chromatic number and the clique number}

Similarly to the independence number and based on the inequality $\chi(G\otimes H)\le \min\{\chi(G),\chi(H)\}$ (see for instance, \cite{ImrKlav00}) we get the next trivial lemma.
\begin{lm}
Let $G$ be a graph and let $\Gamma$ be a family of graphs such that $V(F)=V$, for every $F\in \Gamma$. Consider any function $h:E(G)\longrightarrow\Gamma $. Then,
$\chi(G\otimes_h \Gamma)\le \min\{\chi(G),\chi (h(G))\}.$
\end{lm}

The above upper bound is attained as shown the following example.

\begin{ex}
Let $V=\{x,y,z,t\}$. Consider the graphs $F_i$ on $V$, $i=1,2$, defined by, $E(F_1)=\{xz,yz,zt\}$, and $E(F_2)=\{xy,xz,zt\}$. Let $h: E(K_4)\rightarrow \{F_i\}_{i=1}^2$ be any function that assigns $F_1$ to all of its edges except to one that receives $F_2$. Then, since the graph $K_4\otimes_h\{F_i\}_{i=1}^2$ contains an odd cycle, we have that $3=\chi(K_4\otimes_h\{F_i\}_{i=1}^2)=\chi (h(G))$.
\end{ex}

However, it is not difficult to find examples in which the above upper bound it is not attained.

\begin{ex}
Let $V=\{x,y,z,t\}$ and let $V(K_4)=\{a,b,c,d\}$. Consider the graphs $F_i$ on $V$, $i=1,2$, defined by, $E(F_1)=\{xy,yz,zt\}$, and $E(F_2)=\{xz,xt,yt\}$. Let $h: E(K_4)\rightarrow \{F_i\}_{i=1}^2$ be the function defined by $h(e)=F_1$ if $e\neq ac$ and $h(ac)=F_2$. Consider the function $f:V(K_4\otimes_h\{F_i\}_{i=1}^2)\rightarrow \{0,1,2\}$ defined by
\begin{eqnarray*}
f(a,x)&=& f(b,x)=f(b,z)=f(c,x)=f(d,z)=0, \\
f(a,y)&=& f(b,y)=f(b,t)=f(c,y)=f(d,y)=f(d,t)=1,\\
f(a,z)&=&f(a,t)=f(c,z)=f(c,t)=f(d,x)=2.
\end{eqnarray*}
We have that $h(G)\cong K_4$ and since the graph $K_4\otimes_h\{F_i\}_{i=1}^2$ contains an odd cycle (for instance, the subgraph generated by $\{(a,z), (b,y), (c,x)\}$, we obtain that $3=\chi(K_4\otimes_h\{F_i\}_{i=1}^2)<\chi (h(G))$.
\end{ex}

Related to the clique number, we have the following results.

\begin{lm}\label{lemma_upper_bound_clique_number}
Let $G$ be a graph and let $\Gamma$ be a family of graphs such that $V(F)=V$, for every $F\in \Gamma$. Consider any function $h:E(G)\longrightarrow\Gamma $. Then,
$\omega(G\otimes_h \Gamma)\le \min\{\omega(G),\omega (h(G))\}.$
\end{lm}
\proof Let $\{(a_i,x_i): \ i=1,2,\ldots, k\}$ be a maximal clique of $G\otimes_h \Gamma$. By definition, we have that, $a_ia_j\in E(G)$ and $x_ix_j\in E(h(a_ia_j))$. Thus, the sets $\{a_i: \ i=1,2,\ldots, k\}$ and $\{x_i: \ i=1,2,\ldots, k\}$ are complete subgraphs in $G$ and $h(G)$, respectively.\qed

Let $\Gamma$ be a family of graphs such that $V(F)=V$, for every $F\in \Gamma$, we denote by $\Sigma \Gamma$ the graph with vertex set $V$ and edge-set $\cup_{F\in \Gamma} E(F)$.

\begin{pro}
Let $G$ be a graph and let $\Gamma$ be a family of graphs such that $V(F)=V$, for every $F\in \Gamma$. Then there exists a function $h:E(G)\longrightarrow\Gamma $ such that
$\omega(G\otimes_h \Gamma)=\min\{\omega(G),\omega(\Sigma \Gamma)\}.$
\end{pro}

\proof By Lemma \ref{lemma_upper_bound_clique_number}, we have that $\omega(G\otimes_h \Gamma)\le\min\{\omega(G),\omega(h(G))\}$ and thus, since $E(h(G))\subset E(\Sigma \Gamma)$, $\omega(G\otimes_h \Gamma)\le\min\{\omega(G),\omega(\Sigma \Gamma)\}$. What we have to prove is the reverse inequality. Let the sets $\{a_i: \ i=1,2,\ldots, k\}$ and $\{x_i: \ i=1,2,\ldots, k\}$ be the vertices of complete subgraphs in $G$ and $\Sigma \Gamma$, respectively, where $k=\min\{\omega(G),\omega(\Sigma \Gamma)\}$. By definition, for each pair $i,j$ with $1\le i\le j\le k$ there exists $F_{ij}\in \Gamma$ such that $x_ix_j\in E(F_{ij})$. Then, the function $h:E(G)\longrightarrow\Gamma $
 defined by $h(a_ia_j)=F_{ij}$ produces a complete subgraph with vertices $\{(a_i,x_i):\ i=1,2,\ldots, k\}$ in $G\otimes_h \Gamma$, and thus $\omega(G\otimes_h \Gamma)\ge k$. This proves the result.\qed
\begin{co}
Let $\Gamma$ be a family of graphs such that $V(F)=V$, for every $F\in \Gamma$, and let $n=\omega(\Sigma \Gamma)$. Then there exists a function $h:E(K_n)\longrightarrow\Gamma $ such that
$\chi(K_n\otimes_h \Gamma)=n.$
\end{co}
\subsubsection{Chromatic number in $G\circ_h\Gamma$}
Just like in the case of the lexicographic product, we can obtain an upper bound, not difficult to prove.
\begin{lm}
Let $G$ be a graph of order $n\ge 2$ and $\Gamma$ be a family of graphs. Consider any function $h:V(G)\longrightarrow\Gamma $. Then
$\chi(G\circ_h \Gamma)\le \chi (G)\max_{v\in V(G)}\chi(h(v)).$
\end{lm}
\proof Let $\chi(G)=r$ and $\chi(h(a))=s_a$, $a\in V(G)$. Let $g$ be an $r$-coloring of $G$ and $h_a$ be an $s_a$-coloring of $h(a)$. We claim that $f(a,x)=(g(a),h_a(x))$ defines a coloring of $G\circ_h\Gamma$ with at most $r\max_{v\in V(G)}s_v$ colors. Indeed, suppose that $(a,x)(b,y)\in E(G\circ_h\Gamma)$, then, either $a=b$ and $xy\in E(h(a))$, which implies, since $h_a$ is an $s_a$-coloring, that $h_a(x)\neq h_a(y)$; or, $ab\in E(G)$, but then, since $g$ is an $r$-coloring of $G$, we have $g(a)\neq g(b)$ and the result follows. \qed

The next examples show that the above upper bound is sharp and also, that there exist families of graphs for which the difference between the exact value and the upper bound can be arbitrarily large.

\begin{ex}
Let $C_3=(\{a,b,c\},\{ab,bc,ac\})$ and consider the function $h: V(C_3)\rightarrow \{K_2, K_2\cup K_1\}$ defined by $h(a)=h(b)=K_2$ and $h(c)=K_2\cup K_1$. Then $\chi(C_3\circ_h\{K_2, K_2\cup K_1\})=\chi (C_3)\chi (K_2)$.

Let $C_5=(\{a,b,c,d,e\},\{ab,bc,cd,de,ae\})$ and consider the function $h: V(C_5)\rightarrow \{K_{n},K_2, 2K_1\}$ defined by $h(a)=K_{n}$, $h(b)=K_2$ and $h(c)=h(d)=h(e)=2K_1$. Then $\chi (G)\max_{v\in V(G)}\chi(h(v))-\chi(G\circ_h \Gamma)=2(n-1)$. Indeed, let $\{0,1,\ldots,n-1\}$ and $\{n,n+1\}$ be the colors assigned to $\{(a,x): x\in V(K_n)\}$ and $\{(b,y): y\in V(K_2)\}$ respectively. By assigning $1,n$ and $n+1$ to $\{c\}\times V(2K_1)$, $\{d\}\times V(2K_1)$ and $\{e\}\times V(2K_1)$, respectively, we obtain a $(n+2)$-coloring of $G\circ_h \Gamma$.

\end{ex}

One of the main results found in the study of the chromatic number of the lexicographic product of graphs is the following one due to Geller and Sahl \cite{GelSah75}.
\begin{tm}\cite{GelSah75}\label{chromatic_lexico}
If $\chi (H)=n$ then $\chi(G\circ H)=\chi (G\circ K_n)$.
\end{tm}

In the next lines we generalize Theorem \ref{chromatic_lexico} to the $\circ_h$-product of graphs using similar ideas as the ones found in \cite{GelSah75}. We first recall Proposition 1.20 in \cite{ImrKlav00}.
\begin{pro}\cite{ImrKlav00}\label{char_homo_chromatic_number}
Let $G$ be a graph. Then
 $\chi(G)$ is the smallest integer $n$ for which there exists a homomorphism $G\rightarrow K_n$.
 Moreover, if there exists a homomorphism $G\rightarrow H$, then $\chi(G)\le \chi (H)$.
\end{pro}

\begin{tm}
Let $G$ be a nontrivial graph, $\Gamma$ be a family of graphs and $\{K_{m}\}_{m\in \mathbb{N}}$ be the family of complete graphs. Consider any function $h:V(G)\longrightarrow\Gamma $. Then
$$\chi(G\circ_h\Gamma)=\chi(G\circ_{h'}\{K_{m}\}_{m\in \mathbb{N}}),$$
where $h':V(G)\longrightarrow\{K_{m}\}_{m\in \mathbb{N}}$ is the function defined by $h'(v)=K_n$ if $\chi(h(v))=n$, for every $v\in V(G)$.
\end{tm}

\proof Let $v\in V(G)$ with $\chi(h(v))=n$. By Proposition \ref{char_homo_chromatic_number} there exists a homomorphism $f_v:h(v)\rightarrow K_n$. Thus, we can construct a homomorphism $f$ from  $G\circ_h\Gamma$ onto $G\circ_{h'}\{K_{m}\}_{m\in \mathbb{N}}$ defined by $f(a,x)=(a,f_a(x))$, for each $(a,x)\in V(G\circ_h\Gamma).$ Hence, again by Proposition \ref{char_homo_chromatic_number}, we have that $\chi(G\circ_h\Gamma)\le \chi(G\circ_{h'}\{K_{m}\}_{m\in \mathbb{N}})$.

Conversely, let $f$ be an $r$-coloring of $G\circ_h\Gamma$, with $r=\chi (G\circ_h\Gamma)$. Let $a\in V(G)$, the restriction to the set $\{a\}\times V(h(a))$ contains at least $n=\chi (h(a))$ colors. Choose $n$ of them and a representative vertex in each color class. By connecting each pair of chosen vertices (in case they are not connected), eliminating the extra vertices and repeating the same process for every $a\in V(G)$,  we obtain an $r$-coloring of a graph which is isomorphic to $G\circ_h'\{K_{m}\}_{m\in \mathbb{N}}$, where $h'(a)=K_n$, if $\chi(h(a))=n$. By definition of the chromatic number, we obtain $\chi (G\circ_h'\{K_{m}\}_{m\in \mathbb{N}})\le \chi (G\circ_h\Gamma)$.\qed

The reformulations that have been studied with respect to the chromatic number of the lexicographic product (see \cite{ImrKlav00}) suggest new lines for future research. Suppose that we have a graph $G$, a function $h:V(G)\rightarrow \mathbb{Z^+}$, and we assigne $h(a)$ different colors from the set $\{0,1,2,\ldots,s-1\}$ to each vertex $a$ of $G$ and adjacent vertices receive disjoint sets of colors. In that case, we say that the assignment is a $h$-{\it tuple coloring}. The $h$-{\it chromatic number} $\chi_{h}(G)$ of $G$ is the smallest $s$ such that there is a $h$-{\it tuple coloring} with $s$ colors. When $h$ is constant and equal to $n$, then $h$-{\it tuple coloring} (and the $h$-{\it chromatic number} ) correspond to the $n$-tuple coloring (and the $n^{th}$ chromatic number) that was introduced by Stahl in \cite{Sta76}.

Notice that, similar to what happens with the $nth$ chromatic number, we have that $\chi_{h}(G)=\chi (G\circ_{h'}\{K_{m}\}_{m\in \mathbb{N}})$, where $h'(v)=K_n$ if $h(v)=n$. And we also can establish a relation between $\chi (G\circ_{h'}\{K_{m}\}_{m\in \mathbb{N}})$ and Kneser graphs from a system of sets that we introduce in the following lines.

Let $\{r_i\}_{i\in I}$ be a sequence of positive integers. Denote by $K(\{r_i\}_{i\in I},s)$ the graph that has as vertex set the $r_i$-subsets of a $s$-subset, for each $i\in I$, and two vertices are adjacent if and only if the subsets are disjoint. Clearly, each coloring $c$ of  $G\circ_{h'}\{K_{m}\}_{m\in \mathbb{N}}$ induces a homomorphism $f$ from $G$ onto $K(\{h(a)\}_{a\in V(G)},s)$, defined by $f(a)=\{c(a,x): \ x\in V(h'(a))\}$, where $s$ is the number of colors used and $h(a)=|V(h'(a))|$. Moreover, for every homomorphism $f:V(G)\rightarrow V(K(\{r_i\}_{i\in I},s))$ we obtain an $s$ coloring of $G\circ_{h'}\{K_{m}\}_{m\in \mathbb{N}}$, where $h'(v)=K_n$ if $|f(v)|=n$. Thus, we get the next proposition.

\begin{pro} Let $h: V(G)\rightarrow \mathbb{N}$ be any function and $h':V(G)\rightarrow \{K_{m}\}_{m\in \mathbb{N}}$ be the function defined by $h'(v)=K_{h(v)}$, for all $v\in V(G)$. Then,
$\chi (G\circ_{h'}\{K_{m}\}_{m\in \mathbb{N}})$ is the smallest integer $s$ such that there exists a homomorphism $f$ from $G$ onto $K(\{h(a)\}_{a\in V(G)},s)$ such that $|f(v)|=h(v)$, for all $v\in V(G)$.
\end{pro}
%There is a lower bound for the chromatic number of Kneser graphs from a finite system of sets, that was stated by Dol'nikov's, in terms of the $2$-{\it colorability defect}.

%\begin{open}
%Which one is the chromatic number of $K(\{r_i\}_{i\in I}, s)$?
%\end{open}

\section{Some structural properties}
The next results can be though as some type of associative property for the two products, $\otimes_h$ and $\circ_h$.

\begin{lm}
Let $G$ and $H$ be graphs and let $\Gamma$ be a family of graphs such that  $V(F)=V$ for every $F\in \Gamma$. Then,
\begin{enumerate}
  \item[(i)] For all $h:E(H)\rightarrow \Gamma$ there exists $h':E(G\otimes H)\rightarrow \Gamma$ such that
  $G\otimes (H\otimes_h\Gamma)\cong (G\otimes H)\otimes_{h'}\Gamma.$
  \item[(ii)] For all $h:E(G\otimes H)\rightarrow \Gamma$ with $h((\alpha,a)(\beta,b))=h((\alpha,b)(\beta,a))$, there exist a family $\Gamma'$, with $V(F)=V(H)\times V$, for all $F\in \Gamma'$, and a function $h':E(G)\rightarrow \Gamma'$ such that
  $(G\otimes H)\otimes_h\Gamma\cong G\otimes_{h'}\Gamma'.$
\end{enumerate}
\end{lm}
\proof (i) Let $h:E(H)\rightarrow \Gamma$ and let $h':E(G\otimes H)\rightarrow \Gamma$ be the function defined by $h'((\alpha,a)(\beta,b))=h(ab)$. Then, we have that $V(G\otimes (H\otimes_h\Gamma))=V((G\otimes H)\otimes_{h'}\Gamma)$ and the identity function between the sets of vertices defines an isomorphism of graphs. Indeed, $((\alpha,a),x)((\beta,b),y)\in E((G\otimes H)\otimes_{h'}\Gamma)$ if and only if
$$\left\{\begin{array}{l}
 (\alpha,a)(\beta,b)\in E(G\otimes H)\\
 xy\in E(h'( (\alpha,a)(\beta,b)))
\end{array}\right.\Leftrightarrow\left\{\begin{array}{l}
 \alpha\beta\in E(G)\ \hbox{and} \ ab\in E(H)\\
 xy\in E(h(a,b))\end{array}\right.$$
that is, if and only if $(\alpha,(a,x))(\beta,(b,y))\in E(G\otimes (H\otimes_h\Gamma))$.

(ii) Let $h:E(G\otimes H)\rightarrow \Gamma$ and let $\Gamma'=\{H\otimes_{ h_{\alpha\beta}}\Gamma\}_{\alpha\beta\in E(G)}$, where $h_{\alpha\beta}: E(H)\rightarrow \Gamma$ is the function defined by $h_{\alpha\beta}(ab)=h((\alpha,a)(\beta,b))$. Consider now, the function $h':E(G)\rightarrow \Gamma'$ defined by, $h'(\alpha\beta)=H\otimes_{ h_{\alpha\beta}}\Gamma$. Then, an easy check shows that $V((G\otimes H)\otimes_h\Gamma)=V(G\otimes_{h'}\Gamma')$ and the identity function between the sets of vertices defines an isomorphism of graphs.
\qed

\begin{lm}
Let $G$ and $H$ be graphs and let $\Gamma$ be a family of graphs. Then,
\begin{enumerate}
  \item[(i)] For all $h:V(H)\rightarrow \Gamma$ there exists $h':V(G\circ H)\rightarrow \Gamma$ such that
  $G\circ (H\circ_h\Gamma)\cong (G\circ H)\circ_{h'}\Gamma.$
  \item[(ii)] For all $h:V(G\circ H)\rightarrow \Gamma$ there exists a family $\Gamma'$ and a function $h':E(G)\rightarrow \Gamma'$ such that
  $(G\circ H)\circ_h\Gamma\cong G\circ_{h'}\Gamma'.$
\end{enumerate}
\end{lm}
\proof   (i) Let $h:V(H)\rightarrow \Gamma$ and let $h':V(G\otimes H)\rightarrow \Gamma$ be the function defined by $h'(\alpha,a)=h(a)$. Then, an easy check shows that $V(G\circ (H\circ_h\Gamma))=V((G\circ H)\circ_{h'}\Gamma)$ and the identity function between the sets of vertices defines an isomorphism of graphs. Indeed, $((\alpha,a),x)((\beta,b),y)\in E((G\circ H)\circ_{h'}\Gamma)$ if and only if
$$\left\{\begin{array}{l}
 (\alpha,a)(\beta,b)\in E(G\circ H),\  \hbox{or}\\
(\alpha,a)=(\beta,b) \ \hbox{and}\ xy\in E(h' (\alpha,a)).
\end{array}\right.\Leftrightarrow\left\{\begin{array}{l}
 \alpha\beta\in E(G),\ \hbox{or} \\ \alpha=\beta \ \hbox{and}\  ab\in E(H),\  \hbox{or}\\
\alpha=\beta, \ a=b \ \hbox{and}\ xy\in E(h(a)).\end{array}\right.$$
That is, if and only if $(\alpha,(a,x))(\beta,(b,y))\in E(G\circ (H\circ_h\Gamma))$.

(ii) Let $h:V(G\otimes H)\rightarrow \Gamma$ and let $\Gamma'=\{H\circ_{ h_{\alpha}}\Gamma\}_{\alpha\in V(G)}$, where $h_{\alpha}: V(H)\rightarrow \Gamma$ is the function defined by $h_{\alpha}(a)=h(\alpha,a)$. Consider now, the function $h':V(G)\rightarrow \Gamma'$ defined by, $h'(\alpha)=H\circ_{ h_{\alpha}}\Gamma$. Then, an easy check shows that $V((G\circ H)\circ_h\Gamma)=V(G\circ_{h'}\Gamma')$ and the identity function between the sets of vertices defines an isomorphism of graphs.
\qed

\subsection{On the $\otimes_h$-decomposition for graphs}
Notice that, each graph $G$ admits a trivial decomposition in terms of the $\otimes_h$-product, namely $G\cong L\otimes G$, where $L$ denotes the graph with $|V(L)|=|E(L)|=1$.
We say that $G$ has a nontrivial decomposition with respect the $\otimes_h$-product if there exist a graph $H$ or order at least $2$ (maybe with loops), a family of graphs $\Gamma$ (maybe with loops), with $V(F)=V$ for every $F\in \Gamma$ and a function $h: E(H)\rightarrow \Gamma$ such that $G\cong H\otimes_h\Gamma$. The next result gives necessary and sufficient conditions for the existence of nontrivial $\otimes_h$-decomposition for graphs.

\begin{tm}
Let $G$ be a graph. Then, $G$ has a nontrivial decomposition with respect the $\otimes_h$-product if and only if there exists a partition $V(G)=V_1\cup V_2\ldots \cup V_k$, $k\ge 2$, such that, for each $i,j$ with $1\le i\le j\le k$,
$|V_i|=|V_j|$ and,
 there exist bijective functions $\varphi_{i}: V_1\rightarrow V_i$, such that, for each pair $u,v\in V_1$, we have that
 \begin{equation}\label{condition_decomposition}
   \varphi_{i}(u)\varphi_j(v)\in E(G) \Leftrightarrow \varphi_{i}(v)\varphi_j(u)\in E(G).
 \end{equation}
\end{tm}

\proof Assume that there exist a nontrivial graph $H$, a family of graphs $\Gamma$, with $V(F)=V$ for every $F\in \Gamma$ and a function $h: E(H)\rightarrow \Gamma$ such that $G\cong H\otimes_h\Gamma$. Clearly the $V$-fibers of $H\otimes_h\Gamma$ form a partition of $V(G)$, namely $\cup_{a\in V(H)}\ _aV$. Let $a\in V(H)$. For any vertex $b$ of $H$, consider the function $\varphi_{b}$ defined by $\varphi_{b}(a,x)=(b,x)$. Then, by definition of the $\otimes_h$-product, we have that $(b,x)(c,y)\in E(H\otimes_h\Gamma)$ if and only if  $(b,y)(c,x)\in E(H\otimes_h\Gamma)$. Thus, condition (\ref{condition_decomposition}) holds.

Let us see the sufficiency. Assume that there exists a partition $V(G)=V_1\cup V_2\ldots \cup V_k$, $k\ge 2$, such that, for each $i,j$ with $1\le i\le j\le k$,
$|V_i|=|V_j|$ and,
 there exist bijective functions $\varphi_{i}: V_1\rightarrow V_i$, such that, for each pair $u,v\in V_1$, we have that
  $\varphi_{i}(u)\varphi_j(v)\in E(G)$ if and only if $\varphi_{i}(v)\varphi_j(u)\in E(G)$.
 % Without loss of generality, we can assume that $\varphi_1$ is the identity function.
   Let $V_1=\{x_s\}_{s=1}^l$ and let $H$ be the graph with vertex set $V(H)=\{a_1,a_2,\ldots, a_k\}$ and $a_ia_j\in E(H)$ if and only if $N_G(V_i)\cap V_j\neq \emptyset$, where $
N_G(V_i)=\cup_{v\in V_i} N_G(v)$. For every $a_ia_j\in E(H)$, we consider the graph $F_{ij}$ with vertex set $V_1$ and edge set defined by $x_{s}x_t\in E(F_{ij})$ if and only if $\varphi_i(x_s)\varphi_{j}(x_t)\in E(G)$. Condition (\ref{condition_decomposition}) guarantees that the graph $F_{ij}$ is well defined.
Then, the bijective function $f:V(G)\rightarrow V(H)\times V_1$ defined by $f(v)=(a_i,\varphi_i^{-1}(v))$ if $v\in V_i$, establishes an isomorphism between $G$ and $H\otimes_h \Gamma$, where $\Gamma=\{F_{ij}\}_{a_ia_j\in E(H)}$ and $h:E(H)\rightarrow \Gamma$ is the function defined by $h(a_ia_j)=F_{ij}$.
\qed

Notice that, if we require $H$ to be a loopless graph then we can obtain a similar characterization only by adding the restriction on $V_i$ that says that $V_i$ is formed by independent vertices, for every $i\in \{1,2,\ldots, k\}$. Moreover, if we also require that the family $\Gamma$ does not contain graphs with loops, then we should add the restriction $\varphi_i(u)\varphi_j(u)\notin E(G)$, for each $u\in V_1$ and for each $i,j$ with $1\le i\le j\le k$.

\subsubsection{Non uniqueness}
The next example shows, as it happens with the direct and the lexicographic products, that we do not have a unique decomposition in terms of the $\otimes_h$-product.

\begin{ex}
Let $V=\{x,y,z,t\}$. Consider the graphs $F_i$ on $V$, $i=1,2$, defined by, $E(F_1)=\{xz,yt\}$ and $E(F_2)=\{xt,yz\}$. Let $h: E(C_3)\rightarrow \{F_i\}_{i=1}^2$ be a function in which $F_1$ is assigned to two edges and $F_2$ to the other edge. Then, $$C_3\otimes_h\{F_i\}_{i=1}^2\cong 2C_6\cong 2K_2\otimes C_3.$$

\end{ex}

\noindent {\bf Acknowledgements} The research conducted in this document by the first author has been supported by the Spanish Research Council under project
MTM2011-28800-C02-01 and  by the Catalan Research Council
under grant 2009SGR1387.


\begin{thebibliography}{xxxxx}
\bibitem{AMR} A. Ahmad, F. A. Muntaner-Batle, M. Rius-Font, On the product $\overrightarrow{C_m} \otimes_h \{{\overrightarrow{C}_n, \overleftarrow{C}_n }\}$ and other related topics, {\it Ars Combin.} in press.
\bibitem{BreKlaRal} B. Bre$\check{s}$ar, S. Klav$\check{z}$ar and D. F. Rall, Dominating direct products of graphs, \emph{Discrete Math.} {\bf 307} (2007), 1636--1642.
%\bibitem{CH} G. Chartrand and L. Lesniak, \emph{Graphs and Digraphs},
%second edition. Wadsworth \& Brooks/Cole Advanced Books and
%Software, Monterey (1986).
%\bibitem {E}H. Enomoto, A. Llad\'{o}, T. Nakamigawa and G. Ringel, Super edge-magic graphs, \emph{ SUT J. Math}. {\bf 34} (1998), 105--109.
%\bibitem {F2} R.M. Figueroa-Centeno, R. Ichishima and F.A. Muntaner-Batle, The place of super edge-magic labelings among other classes of labelings, \emph{Discrete Math.} {\bf 231} (1--3) (2001), 153--168.
%\bibitem{FIM02}R.M. Figueroa-Centeno, R. Ichishima and F.A. Muntaner-Batle, Magical coronations of graphs, \textit{Australasian J. Combin.} \textbf{ 26} (2002) 199--208.
%\bibitem {FIM05b}R. M. Figueroa-Centeno, R. Ichishima and F. A. Muntaner-Batle.\newblock On super edge-magic labelings of certain disjoint unions of graphs. \newblock  \emph{Australas. J. Combin.} {\bf 32} (2005), 225--242.
%\bibitem{F3} R.M Figueroa -Centeno, R. Ichishima and F.A. Muntaner-Batle,  On super edge-magic graphs, \emph{Ars Combin.} {\bf 64} (2002), 81--95.
%\bibitem {FIMO}R.M. Figueroa-Centeno, R. Ichishima, F.A. Muntaner-Batle
%and A. Oshima, A magical approach to some labeling conjectures, \emph{Discuss. Math. Graph Theory} {\bf 311} (2011), 79--113.
\bibitem {F1}R.M. Figueroa-Centeno, R. Ichishima, F.A. Muntaner-Batle
and M. Rius-Font, Labeling generating matrices, \emph{J. Comb. Math. and Comb. Comput.} {\bf 67} (2008), 189--216.
\bibitem{GelSah75} D. Geller and S. Sahl, The chromatic number and other functions of the lexicographic produt, {\it J. Combin. Theory Ser. B}, {\bf 19} (1975), 87--95.
\bibitem{GraKhe95} S. Gravier and A. Khelladi, On the domination number of cross product of two Hamiltonian graphs, {\it Discrete Math.}
{\bf 145} (1995), 273--277.

%\bibitem{HolMcqMcq09} J. Holden, D. McQuillan and J.M. McQuillan, A conjecture on strong magic labelings of $2$-regular graphs, \emph{Discrete Math.} {\bf 309} (2009) 4130--4136.
%\bibitem{HLEK03} P. Holme, F. Liljeros, G.R. Edling and B.J. Kim, Network bipartivity, {\it Phys. Rev. E.} 68, 056107 (2003), 12pp.
\bibitem{ILMR} R. Ichishima, S.C. L\'opez, F. A. Muntaner-Batle, M. Rius-Font, The power of digraph products applied to labelings, {\it Discrete Math.} {\bf 312} (2012) 221-228. %DOI:10.1016/j.disc.2011.08.025
%\bibitem{Imr98} W. Imrich, Factoring cardinal product praphs in polynomial time, {\it Discrete Math.} {\bf 192} (1998), 119--144.
\bibitem{ImrKlav00} W. Imrich and S. Klav$\check{z}$ar, Product Graphs:
Structure and Recognition (John Wiley $\&$ Sons, Inc. New York, 2000).
%
% \bibitem {G}J.A. Gallian, A dynamic survey of graph labeling,\emph{Electron. J. Combin.} \textbf{14} (2011), $\sharp $\textbf{DS6}.
%\bibitem{Gra97} S. Gravier, Hamiltonicity of the cross product of two Hamiltonian graphs, {\it Discrete Math.}
%{\bf 170} (1997), 253--257.
%\bibitem{GodSla98} R.D. Godbold and P. J. Slater, All cycles are edge-magic, \emph{Bull. Inst. Combin Appl.}{\bf 22} (1998), 93--97.

\bibitem{KlaZam} S.  Klav$\check{z}$ar and B. Zmazek, On a Vizing-like conjecture for direct product of graphs, {\it Discrete Math.} {\bf 156} (1996), 243--246.

%\bibitem {K1}A. Kotzig and A. Rosa, Magic valuations of finite graphs,
%\emph{Canad. Math. Bull.} {\bf 13} (1970), 451--461.
\bibitem{LopMunRiu1}  S.C. L\'opez, F. A. Muntaner-Batle, M. Rius-Font, Bi-magic and other generalizations of super edge-magic labelings, \emph{B. Aust. Math. Soc}.\textbf{ 84} (2011) 137--152.
%\bibitem{LopMunRiu1}  S.C. L\'opez, F. A. Muntaner-Batle, M. Rius-Font, Bi-magic and other generalizations of super edge-magic labelings, \emph{B. Aust. Math. Soc}.\textbf{ 84} (2011) 137--152.
\bibitem{LMR_LabCons}  S.C. L\'opez, F. A. Muntaner-Batle, M. Rius-Font, Labeling constructions using the $\otimes_h$-product,  submitted to \emph{Discrete Applied Math.}
\bibitem{Meki10} G. Meki$\check{s}$, Lower bounds for the dominating number and the total domination number of direct product graphs, {\it Discrete Math.} {\bf 310} (2010), 3310-3317.
\bibitem{NowRall96} R. J. Nowakowski and D. F. Rall, Associative graph products and their independence, domination and coloring numbers, {\it Discuss. Math. Graph Theory} {\bf 16} (1996), 53--79.
\bibitem{Sta76} S. Stahl, $n$-tuple colorings and associated graphs, {\it J. Combin. Theory Ser. B} {\bf 20} (1976), 185--203.
%\bibitem {Tar05} C. Tardif, The fractional chromatic number of the categorical product of graphs, {\it Combinatorica}
%{\bf 25} (5) (2005), 625--632.
\bibitem {YangXu} C. Yang and J.-M. Xu, Connectivity of lexicographic product and direct product of graphs,  {\it Ars Combinatoria}, in press.
\bibitem{Wei62} P.W. Weichsel, The Kronecker product of graphs, {\it Proc. Amer. Math. Soc.} {\bf 13} (1962), 47--52.
    \bibitem{WhiRus12} A.N. Whitehead and B. Russell (1912) {\it Principia Mathematica}, volume 2. Cambridge University Press, Cambridge.
\bibitem{Sab61} G. Sabidussi, Graph derivatives, {\it Math. Z.} {\bf 76} (1961) 385--401.
\end{thebibliography}
\end{document}